\documentclass[review]{elsarticle}
\usepackage{lineno,hyperref}
\modulolinenumbers[5]
\journal{}

\usepackage{amsmath,amsthm}
\usepackage{amssymb,latexsym}
\usepackage{graphicx}
\usepackage[labelsep=endash]{caption}
\usepackage{mathptmx}
\usepackage{grffile}
\usepackage{bigints}

\usepackage{subfig}
\newtheorem{thm}{Theorem}[section]
\newtheorem{cor}{Corollary}[section]
\newtheorem{lem}{Lemma}[section]

\theoremstyle{definition}
\newtheorem{dfn}{Definition}[section]
\newtheorem{rem}{Remark}[section]

\begin{document}
\begin{frontmatter}
\title{Generalized least square homotopy perturbations for system of fractional partial differential equations}
\author[label1]{Rakesh Kumar}
\ead[a]{rakesh.lect@gmail.com}
\author[label2]{Reena Koundal}
\address[label1]{Department of Mathematics, Central University of Himachal Pradesh, Dharamshala, India}
\address[label2]{Department of Mathematics, Central University of Himachal Pradesh, Dharamshala, India}
\begin{abstract}
In this paper, generalized aspects of least square homotopy perturbations are explored to treat the system of non-linear fractional partial differential equations and the method is called as generalized least square homotopy perturbations (GLSHP). The concept of partial fractional Wronskian is introduced to detect the linear independence of functions depending on more than one variable through Caputo fractional calculus. General theorem related to Wronskian is also proved. It is found that solutions converge more rapidly through GLSHP in comparison to classical fractional homotopy perturbations. Results of this generalization are validated by taking examples from nonlinear fractional wave equations.  
\end{abstract}

\end{frontmatter}
\section{Introduction}
Fractional calculus is one of the highly specialized branch of applied mathematics which investigates the properties and applications of derivatives and integrals of arbitrary order (real or complex) (\cite{Miller1993}, \cite{Saamko1993}). The definitions of fractional derivative are given by many authors but most popular definitions are due to Gr\"{u}nwald-Letnikov, Riemann-Liouville and Caputo (\cite{Oldham 1974}, \cite{Podlubny1999}). During the last few years, fractional partial differential equations have been the topic of various researchers for accurate modeling and analysis owing to the utility in plasma physics, biology, fluid mechanics and many more \cite{thabetkendre2017}. It is still a challenging task to obtain analytical or approximate solutions to nonlinear fractional partial differential equations. Therefore, mathematicians and scientists are in the search of accurate methods which provide convergent and stable solutions to these equations.\\
Some techniques for the solution of fractional partial differential equations can be seen in literature such as homotopy perturbation method \cite{momani2007}, variational iteration method and decomposition method \cite{odibat2008}, homotopy analysis method \cite{jafari2009} homotopy perturbation technique \cite{sayed2012}, Laplace transform method \cite{jafari2013}, variational iteration with Pade approximations \cite{turut2013}, corrected Fourier series\cite{zainal2014} and fractional complex transformation \cite{zayed2016}. \\
In this paper, we consider a space fractional wave equation in the following form \cite{Ghazanfari2011}:
\begin{align}\label{a.1}
\frac{\partial^{2\alpha}u}{\partial x^{2\alpha}}-u\frac{\partial^{2}u}{\partial t^{2}} = \phi(x,t),\hspace{0.3cm}1<\alpha\leqslant 1,
\end{align}
subject to initial conditions
\begin{align}\label{b.1}
u(0,t)= f_1(t), \hspace{0.3cm} \dfrac{\partial}{\partial x}u(0,t)= g_1(t).
\end{align}
Consider a space fractional system of wave equations as \cite{Biazar2013}:
\begin{align*}
\frac{\partial^{2\alpha}u}{\partial x^{2\alpha}}-v\frac{\partial^{2}u}{\partial t^{2}}-u\frac{\partial^{2}v}{\partial t^{2}} = \phi (x,t),\hspace{0.3cm}0<\alpha\leqslant 1,
\end{align*}
\begin{align}\label{a.2}
\frac{\partial^{2\beta}v}{\partial x^{2\beta}}-v\frac{\partial^{2}v}{\partial t^{2}}+u\frac{\partial^{2}u}{\partial t^{2}} = \psi (x,t),\hspace{0.3cm}0<\beta\leqslant 1.
\end{align}
subject to initial conditions
\begin{align}
u(0,t)= f_1(t), \hspace{0.3cm} \dfrac{\partial}{\partial x}u(0,t)= g_1(t),\\
v(0,t)= f_2(t), \hspace{0.3cm} \dfrac{\partial}{\partial x}v(0,t)= g_2(t).
\end{align} 
In above equations, $\phi(x,t)$, $\psi(x,t)$, $f_i(t)$ and $g_i(t)$ for $i=1,2$ are given functions.
Wei \cite{wei2017} analyzed a local finite difference scheme for the diffusion-wave equation of fractional order. A fourth order compact scheme which preserve energy was devised by Diaz et al. \cite{diaz2018} for the solution of fractional nonlinear wave equations. Weak solutions for time fractional diffusion equations were examined by Yamamoto \cite{yamamoto2018}. Other schemes which prescribe the information regarding the solution of fractional wave equations can be studied from \cite{Khaled2005}-\cite{deng2015}.\\
In previous years, scientists and applied mathematicians have been attracted towards modifications of classical homotopy perturbations in order to achieve accelerated accuracy. Some modifications can be seen in \cite{Wu2010} and \cite{Ullah2015}. A variational form for homotopy perturbation iteration method was developed for fractional diffusion equation by Guo et al. \cite{guo2013}. Recently, Constantin and Caruntu \cite{BOTA2017} presented least square homotopy perturbations for non-linear ordinary differential equations. But this technique was not suitable to handle fractional partial differential equations as convergent solutions are not expected. Therefore, our main objective is to generalize this idea of the coupling of fractional homotopy perturbations and least square approximations, and to propose fractional partial Wronskian.
The paper is organized in the following sequence: In section 2, basic definitions of fractional calculus are presented. The definition and theory of fractional partial Wronskian is also developed. In section 3, basic theory of generalized least square homotopy perturbations are proposed together with necessary definitions and lemma. Section 4 deals with numerical examples. Lastly conclusion is derived.
\section{Fundamentals of fractional calculus}
We provide some basic definitions, properties and theorem related to fractional calculus which will be used at later stage (\cite{Miller1993}-\cite{Podlubny1999}).
\begin{dfn}
A real function $f(t)$, $t>0$ is said to be in the space $ C_{\mu} $, $ \mu \in \mathbb{R}$ if there exists a real number $ p>\mu $, such that $ f(t)=t^{p}f_{1}(t) $ where $ f_{1}(t) \in C(0,\infty)$, and is said to be in the space $ C_{\mu}^{n} $ if and only if $ f^{n}\in C_{\mu} $, $ n\in \mathbb{N} $ .
\end{dfn}
\begin{dfn}
The (left sided) Riemann-Liouville fractional integral operator of order $ \alpha \geqslant 0 $ for a function $ f\in C_{\mu} $, $ \mu \geqslant -1 $ is defined as:
\begin{align*}
J^{\alpha} f(t)=
\begin{cases}
\dfrac{1}{\Gamma (\alpha)} \bigintsss_{0}^{t}(t-\tau)^{\alpha-1}f({\tau}) d \tau,&  \alpha > 0, \hspace{0.2cm} t>0,\\
 f(t), \hspace{0.2cm}  \alpha = 0.\\                        
\end{cases}
\end{align*}
\end{dfn}
\begin{dfn}
The (left sided) Caputo fractional derivative of $ f $, $ f\in C_{-1}^{n}$, $ n\in \mathbb{N} \bigcup \lbrace 0 \rbrace $ is defined as:
\begin{align*}
 D_*^{\alpha} f(t)=
\begin{cases}
I^{n-\alpha}f^{n}(t),&  n-1<\alpha < n, \hspace{0.2cm} n \in N,\\
   \dfrac{d^{n}}{d t^{n}}f(t),& \alpha =n.                    
\end{cases}
\end{align*}
\end{dfn}
Some basic properties of $ J^{\alpha}$ and $D_*^\alpha$ are mentioned below for function $ f\in C_{\mu} $ with $ \mu \geq -1$, $ \alpha,\beta \geq 0$ and $ \gamma \geq 1 $:
\begin{align*}
&i)\hspace{0.2cm} J^{\alpha}J^{\beta} f(t)=J^{\alpha+\beta} f(t), \hspace{0.4cm} ii) \hspace{0.2cm} J^{\alpha}J^{\beta}f(t) =J^{\beta} J^{\alpha}f(t),& \hspace{0.4cm} iii) \hspace{0.2cm} J^{\alpha} t^{\gamma}= \dfrac{\Gamma (\gamma+1)}{\Gamma(\alpha+\gamma+1)}t^{\alpha+\gamma},\\
&iv)\hspace{0.2cm} D_*^{\alpha}J^{\alpha}f(t)=f(t), \hspace{0.4cm} v) \hspace{0.2cm}  J^{\alpha}D_*^{\alpha}f(t) = f(t)-\sum_{i=0}^{n-1} f^{i}(0^{+})\dfrac{t^{i}}{i!},&\\
&vi)\hspace{0.2cm} D_*^{\alpha} t^{\gamma}=\begin{cases}
0,&  \gamma\leq\alpha-1, \\
   \dfrac{\Gamma (\gamma+1)}{\Gamma(\gamma-\alpha+1)}t^{\gamma-\alpha},& \gamma>\alpha-1.                   
\end{cases}&
\end{align*}
\begin{dfn}
 For $n$ to be the smallest integer that exceeds $ \alpha $, the Caputo fractional partial derivative of order $ \alpha>0 $ with respect to $x$ is defined as:
\begin{align*}
 \dfrac{\partial^{\alpha}u(x,t)}{\partial x^{\alpha}}=
\begin{cases}
\dfrac{1}{\Gamma (n-\alpha)} \bigintsss_{0}^{x}(x-\tau)^{\alpha-1}\dfrac{\partial^{n}u(\tau,t)}{\partial x^{n}} d \tau,&   n-1<\alpha < n, \\
   \dfrac{\partial^{n}u(x,t)}{\partial x^{n}},&  \alpha=n\in \mathbb{N}.                   
\end{cases}
\end{align*}
\end{dfn}
\begin{dfn} \label{2.5}
Let $\Phi_{1}, \Phi_{2},...,\Phi_{n}$ be $n$ functions of variables $x$ and $t$ which are defined on domain $\Omega$, then fractional partial Wronskian of $\Phi_{1}, \Phi_{2},...,\Phi_{n}$ is 
\[
W^{\alpha}[\Phi_{1}, \Phi_{2},\dots,\Phi_{n}] =
\begin{vmatrix}
\Phi_{1}& \Phi_{2} &\Phi_{3} & \dots & \Phi_{n} \\ 
D^{\alpha} \Phi_{1} & D^{\alpha}\Phi_{2} & D^{\alpha} \Phi_{3} & \dots & D^{\alpha} \Phi_{n} \\
D^{2\alpha} \Phi_{1} & D^{2\alpha}\Phi_{2} & D^{2\alpha} \Phi_{3} & \dots & D^{2\alpha}\Phi_{n} \\
\vdots& \vdots& \vdots& \vdots& \vdots\\
D^{\alpha(n-1)} \Phi_{1} & D^{\alpha(n-1)} \Phi_{2} &  D^{\alpha(n-1)} \Phi_{3}&  \dots & D^{\alpha(n-1)} \Phi_{n}
\end{vmatrix}
,\]
where $ D^{\alpha}(\Phi_{i})=\left(\dfrac{\partial}{\partial t} + \dfrac{\partial^{\alpha}}{\partial x^{\alpha}}\right) (\Phi_{i}) $ and $D^{n \alpha}=D^{\alpha}.D^{\alpha}...D^{\alpha} (n-\text{times})$ for $ 0<\alpha\leq 1$ and $i=1,2,3........n. $
\end{dfn}

\begin{thm}
If the fractional partial Wronskian of $ n $ functions $ \Phi_{1}(x,t) $, $ \Phi_{2}(x,t) $, $\dots$, $ \Phi_{n}(x,t) $ is non zero, at least at one point of the domain $\Omega=[a,b]\times [a,b]$, then functions $ \Phi_{1}(x,t) $, $ \Phi_{2}(x,t) $,...,$ \Phi_{n}(x,t) $ are said to be linearly independent. 
\end{thm}
\begin{proof}
Assume that the functions $ \Phi_{1}(x,t) $, $ \Phi_{2}(x,t) $, $\dots$ ,$ \Phi_{n}(x,t) $ are linearly dependent on the interval $ [a,b] $. So, there exists constants (at least one of them is nonzero) $ k_{1} $, $ k_{2} $, $\dots$ , $ k_{n} $ such that
\begin{align}\label{w.2}
k_{1} \Phi_{1}(x,t)+k_{2} \Phi_{2}(x,t)+ \dots+k_{n}\Phi_{n}(x,t)=0, \hspace{0.3cm} \forall  (x,t) \in [a,b]\times [a,b].
\end{align}
Operating above equation repeatedly $(n-1)$-times by the operator  $ D^{\alpha}=\dfrac{\partial}{\partial t} + \dfrac{\partial^{\alpha}}{\partial x^{\alpha}}$, we obtain the following set of equations:
\begin{align}\label{t1}
k_{1}D^{\alpha} \Phi_{1}(x,t)+k_{2}D^{\alpha} \Phi_{2}(x,t)+& \dots+k_{n}D^{\alpha}\Phi_{n}(x,t)=0, \hspace{0.3cm} \forall  (x,t) \in [a,b]\times [a,b],\\
k_{1}D^{2\alpha} \Phi_{1}(x,t)+k_{2}D^{2\alpha} \Phi_{2}(x,t)+& \dots+k_{n}D^{2\alpha}\Phi_{n}(x,t)=0, \hspace{0.3cm} \forall  (x,t) \in [a,b]\times [a,b],\\
\vdots \nonumber \\
k_{1}D^{(n-1)\alpha} \Phi_{1}(x,t)+k_{2}D^{(n-1)\alpha} \Phi_{2}(x,t)+& \dots+k_{n}D^{(n-1)\alpha}\Phi_{n}(x,t)=0, \hspace{0.3cm} \forall  (x,t) \in [a,b]\times [a,b]. \label{t2}
\end{align}
Let $ (x_{0},t_{0}) \in [a,b]\times[a,b] $ be an arbitrary point such that
\begin{align*}
W^{\alpha}[\Phi_{1}, \Phi_{2},\dots,\Phi_{n}](x_{0},t_{0}) =
\begin{vmatrix}
\Phi_{1}(x_{0},t_{0})& \Phi_{2}(x_{0},t_{0}) &\Phi_{3}(x_{0},t_{0}) & \dots & \Phi_{n}(x_{0},t_{0}) \\ 
D^{\alpha} \Phi_{1}(x_{0},t_{0}) & D^{\alpha} \Phi_{2}(x_{0},t_{0}) & D^{\alpha} \Phi_{3}(x_{0},t_{0}) & \dots & D^{\alpha} \Phi_{n}(x_{0},t_{0}) \\
D^{2\alpha} \Phi_{1}(x_{0},t_{0}) & D^{2\alpha} \Phi_{2}(x_{0},t_{0}) & D^{2\alpha} \Phi_{3}(x_{0},t_{0}) & \dots & D^{2\alpha} \Phi_{n}(x_{0},t_{0}) \\
\vdots& \vdots & \vdots & \vdots & \vdots \\
D^{\alpha(n-1)} \Phi_{1}(x_{0},t_{0}) & D^{\alpha(n-1)} \Phi_{2} (x_{0},t_{0})&  D^{\alpha(n-1)} \Phi_{3}(x_{0},t_{0})&  \dots & D^{\alpha(n-1)} \Phi_{n}(x_{0},t_{0})
\end{vmatrix}
\neq 0.
\label{t}
\end{align*}
Equations (\ref{w.2})-(\ref{t2}) at the point $ (x_{0},t_{0})$ have the following representations:
\begin{equation} \label{eq1}
\begin{split}
k_{1} \Phi_{1}(x_{0},t_{0})+k_{2} \Phi_{2}(x_{0},t_{0})&+\dots+k_{n} \Phi_{n}(x_{0},t_{0})=0, \\
k_{1}D^{\alpha} \Phi_{1}(x_{0},t_{0})+k_{2}D^{\alpha}\Phi_{2}&(x_{0},t_{0})+\dots+k_{n}D^{\alpha} \Phi_{n}(x_{0},t_{0})=0,\\
&\vdots\\
k_{1}D^{(n-1)\alpha} \Phi_{1}(x_{0},t_{0})+k_{2}D^{(n-1)\alpha} \Phi_{2}&(x_{0},t_{0})+\dots+k_{n}D^{(n-1)\alpha}\Phi_{n}(x_{0},t_{0})=0.
\end{split}
\end{equation}
The determinant of the coefficient matrix of system (\ref{eq1}) is the Wronskian $W^{\alpha}[\Phi_{1}, \Phi_{2},\dots,\Phi_{n}]$ at the point $(x_{0},t_{0})$. Also $W^{\alpha}[\Phi_{1}, \Phi_{2},\dots,\Phi_{n}](x_{0},t_{0})\neq 0 $, therefore, $k_{1}=k_{2}=\dots=  k_{n}=0.$ Which contradicts our hypothesis.\\
Hence, the functions $ \Phi_{1}(x,t) $, $ \Phi_{2}(x,t) $, $ \dots, $ $ \Phi_{n}(x,t) $ are linearly independent. 
\end{proof}
\section{General procedure  of generalized least square homotopy perturbations (GLSHP)}
We consider the following system of fractional partial differential equations along with initial and boundary conditions as: 
\begin{align}\label{1.1}
L_{1}^{\alpha}(u(x,t))+\mathcal{N}(u(x,t),v(x,t))-f(x,t)=0,\hspace{0.3cm} x,t\in\Omega = [0,1]\times[0,1], \hspace{0.3cm} 0<\alpha\leqslant 1,\\ 
L_{2}^{\beta}(v(x,t))+\mathcal{M}(u(x,t),v(x,t))-g(x,t)=0, \hspace{0.3cm} x,t\in\Omega = [0,1]\times[0,1],  \hspace{0.3cm} 0<\beta\leqslant 1 \label{1.2}
\end{align}
with
\begin{align}\label{1.1a}
B_{1}\left( u,\frac{\partial u}{\partial n}\right) =0, \hspace{0.4cm}  B_{2}\left( v,\frac{\partial v}{\partial n}\right) = 0, \hspace{0.3cm} x,t\in  \Gamma
\end{align}
In above equations, $L_{1}^{\alpha}$, $ L_{2}^{\beta} $ are linear operators, $\mathcal{N}$, $ \mathcal{M} $ are non-linear operators and $u(x,t)$, $v(x,t)$ are unknown functions. Whereas $B_{1}$, $ B_{2} $ are boundary operators, $f(x,t)$ and $g(x,t)$ are known functions and $ n $ is a normal vector on $ \Gamma $ (boundary of $ \Omega $). \\
Following classical homotopy perturbations, we construct two homotopies as: $U(x,t;p):\Omega \times [0,1]\rightarrow \mathbb{R}$ and $V(x,t;p):\Omega \times [0,1]\rightarrow \mathbb{R}$ which satisfy
\begin{align}\label{1.5}
H_{1}(U(x,t;p),p)=(1-p)\left[L_{1}^{\alpha}(U(x,t;p))-f(x,t)\right]+p\left[L_{1}^{\alpha}(U(x,t;p))+\mathcal{N}(U(x,t;p),(V(x,t;p))-f(x,t)\right] =0,
\end{align}
\begin{align}\label{1.6}
H_{2}(V(x,t;p),p)=(1-p)\left[L_{2}^{\beta}(V(x,t;p))-g(x,t)\right]+p\left[L_{2}^{\beta}(V(x,t;p))+\mathcal{M}(U(x,t;p),(V(x,t;p))-g(x,t)\right]=0.
\end{align}
Putting $p = 0$ in above two equations, we should receive: $ L_{1}^{\alpha}(U(x,t;0))=u_{0}(x,t) $, $  L_{2}^{\beta}(V(x,t;0))=v_{0}(x,t) $ where initial guesses $ u_{0}(x,t) $ and $ v_{0}(x,t) $ will be obtained by solving 
\begin{align}\label{1.na}
L_{1}^{\alpha}(u_{0}(x,t)-f(x,t))=0,
\end{align}
\begin{align}\label{1.naa}
L_{2}^{\beta}(v_{0}(x,t)-g(x,t))=0
\end{align}
with 
\begin{align*}
B_{1}\left( u_{0},\frac{\partial u_{0}}{\partial n}\right) =0, \hspace{0.4cm}  B_{2}\left( v_{0},\frac{\partial v_{0}}{\partial n}\right) = 0.
\end{align*}
The solutions to equations (\ref{1.na}) and (\ref{1.naa}) are 
\begin{align}\label{1.6a}
u_{0}(x,t) = \sum_{i=0}^{n-1}u_{0}^{i}(0^{+})\dfrac{x^{i}}{i!}+J^{\alpha}f(x,t),
\end{align}
\begin{align}\label{1.6aa}
v_{0}(x,t) = \sum_{i=0}^{n-1}v_{0}^{i}(0^{+})\dfrac{x^{i}}{i!}+J^{\beta}g(x,t).
\end{align}
In view of HPM, we use the homotopy parameter $ p $ to write the expansions for $U(x,t,p)$ and $ V(x,t,p) $ as:
\begin{align}\label{1.7}
U(x,t,p)= u_{0}(x,t)+\sum_{r\geqslant 1}u_{r}(x,t)p^{r},
\end{align}
\begin{align}\label{1.8}
V(x,t,p)= v_{0}(x,t)+\sum_{h\geqslant 1}v_{h}(x,t)p^{h}. 
\end{align}
Utilizing equations (\ref{1.7}) and (\ref{1.8}) in equations (\ref{1.5}) and (\ref{1.6}),  and then equating the coefficients of corresponding power of $ p $ on both sides, we have
\begin{align}\label{1.9}
\begin{aligned}
L_{1}^{\alpha}(u_{r}(x,t))=&-\mathcal{N}_{r-1}[\lbrace u_{0}(x,t),v_{0}(x,t)\rbrace,\lbrace u_{1}(x,t)v_{1}(x,t)\rbrace,......,\lbrace u_{r-1}(x,t),v_{h-1}(x,t)\rbrace],
\end{aligned}
\end{align}
\begin{align}\label{1}
L_{2}^{\beta}(v_{h}(x,t))=&-\mathcal{M}_{h-1}[\lbrace u_{0}(x,t),v_{0}(x,t)\rbrace,\lbrace u_{1}(x,t),v_{1}(x,t)\rbrace,,......,\lbrace u_{r-1}(x,t,)v_{h-1}(x,t)\rbrace],
\end{align}
\begin{align*}
B_{1}\left( u_{r},\frac{\partial u_{r}}{\partial n}\right) =0, \hspace{0.4cm}  B_{2}\left( v_{h},\frac{\partial v_{h}}{\partial n}\right) = 0 \hspace{0.3cm}r,h \geqslant 1.
\end{align*}
In above equations, nonlinear operators $\mathcal{N}_{k}(x,t)$ and $ \mathcal{M}_{l}(x,t) $ are coefficients of $ p^{k} $ and $ p^{l} $ in following equations
\begin{align}\label{1.11}
\begin{aligned}
\mathcal{N}(x,t)=& \mathcal{N}_{0} [u_{0}(x,t),v_{0}(x,t)]  +p\mathcal{N}_{1}[ u_{0}(x,t),u_{1}(x,t),v_{0}(x,t),v_{1}(x,t)]  +p^{2}\mathcal{N}_{2}[ u_{0}(x,t),u_{1}(x,t) &\\
&,u_{2}(x,t),v_{0}(x,t),v_{1}(x,t),v_{2}(x,t)]+  \dots,&
\end{aligned}
\end{align}
\begin{align}\label{1.12}
\begin{aligned}
\mathcal{M}(x,t)= &\mathcal{M}_{0}[u_{0}(x,t),v_{0}(x,t)]+p\mathcal{M}_{1}[u_{0}(x,t),u_{1}(x,t),v_{0}(x,t),v_{1}(x,t)]+p^{2}\mathcal{M}_{2}[u_{0}(x,t),u_{1}(x,t)&\\
&,u_{2}(x,t),v_{0}(x,t),v_{1}(x,t),v_{2}(x,t)]+\dots,&
\end{aligned}
\end{align}
where $ r-1=k >0$, $ h-1=l>0 $.\\
Here classical fractional homotopy perturbations will give the $r^\text{(th)}$-order and $h^\text{(th)}$-order approximate solutions (for $p=1$) as
\begin{align*}
u^{r}=u_{0}+u_{1}+u_{2}+\dots+u_{r} \,\,\,\,\text{and}\,\,\,\, v^{h}=v_{0}+v_{1}+v_{2}+\dots+v_{h}.
\end{align*}
But here to move further, we are only interested in exploiting initial iterations of fractional homotopy perturbations $ u_{r} $ and $ v_{h} $, for $ r,h\geq 0$ (generally $0^\text{(th)}$-order or $1^\text{(st)}$-order ),  and rest approximations can be ignored. For this, we pick up linearly independent functions from these initial iterations after verifying that
\[
W^{\alpha}_{1}[u_{0}, u_{1}, u_{2},\dots,u_{r}] =
\begin{vmatrix}
u_{0}& u_{1} &u_{2} & \dots & u_{r} \\ 
D^{\alpha} u_{0} & D^{\alpha} u_{1} & D^{\alpha} u_{2} & \dots & D^{\alpha} u_{r} \\
D^{2\alpha} u_{0} & D^{2\alpha} u_{1} & D^{2\alpha} u_{2} & \dots & D^{2\alpha} u_{r} \\
\vdots& \vdots& \vdots& \vdots& \vdots\\
D^{\alpha(n-1)} u_{0} & D^{\alpha(n-1)} u_{1} &  D^{\alpha(n-1)} u_{2}&  \dots & D^{\alpha(n-1)} u_{r}
\end{vmatrix}
\neq0
\]
and
\[
W^{\beta}_{2}[v_{0}, v_{1}, v_{2},\dots,v_{h}] =
\begin{vmatrix}
v_{0}& v_{1} &v_{2} & \dots & v_{h} \\ 
\mathcal{D^{\beta}} v_{0} & \mathcal{D^{\beta}} v_{1} &\mathcal{D^{\beta}} v_{2} & \dots & \mathcal{D^{\beta}} v_{h} \\
\mathcal{D}^{2\beta} v_{0} & \mathcal{D}^{2\beta} v_{1} & \mathcal{D}^{2\beta} v_{2} & \dots &\mathcal{D}^{2\beta} v_{h} \\
\vdots& \vdots& \vdots& \vdots& \vdots\\
\mathcal{D}^{\beta(n-1)} v_{0} & \mathcal{D}^{\beta(n-1)} v_{1} & \mathcal{D}^{\beta(n-1)} v_{2}&  \dots & \mathcal{D}^{\beta(n-1)} v_{h}
\end{vmatrix}
\neq0.
\]
Here, $ D^{\alpha} $=$\left(\dfrac{\partial}{\partial t}+\dfrac{\partial^{\alpha}}{\partial x^{\alpha}} \right)$ and  $ \mathcal{D^{\beta}} $=$\left(\dfrac{\partial}{\partial t}+\dfrac{\partial^{\beta}}{\partial x^{\beta}} \right)$, where $ \alpha,\beta $ lies between $ 0<\alpha,\beta\leqslant 1. $\\
Now, let $Q_{r}= \lbrace \phi_{0},\phi_{1},\dots,\phi_{r}\rbrace$ and $W_{h} =\lbrace \varphi_{0},\varphi_{1},\dots,\varphi_{h}\rbrace$, where $r,h= 0,1,2\dots$ and $ \phi_{0}=u_{0},\phi_{1}=u_{1},\dots,\phi_{r}=u_{r}$, and $ \varphi_{0}=v_{0},\varphi_{1}=v_{1},\dots,\varphi_{h}=v_{h}$. Here $ Q_{r-1}\subseteq Q_{r} $,  $ W_{h-1}\subseteq W_{h} $ and elements of $ Q_{r} $, $ W_{h} $ are linearly independent in the vector space of continuous functions defined over $ \mathbb{R} $.\\
Let us assume $\tilde{u}=\sum_{i=0}^{r}K_{r}^{i}\phi_{i}$ and $\tilde{v}=\sum_{j=0}^{h}D_{h}^{j}\varphi_{j}$ for $ r,h\geqslant 0$ as the approximate solution to set of equations (\ref{1.1}) and (\ref{1.2}).\\
{\bf{Basic algorithm to calculate the coefficients of approximations:}}\\ 
Now we shall employ \textit{least square} procedure to calculate these coefficients in the approximations for which we require the following definitions and lemma:
\begin{dfn} We define remainders $ R_{1}^{\alpha,\beta} $, $ R_{2} ^{\alpha,\beta}$ for system of partial differential equations (\ref{1.1}), (\ref{1.2}) as:
\begin{align}\label{1.3}
R_{1}^{\alpha,\beta}(t,x,\tilde{u},\tilde{v})=L_{1}^{\alpha}(\tilde{u}(x,t))+\mathcal{N}(\tilde{u}(x,t), \tilde{v}(x,t))-f(x,t) , \hspace{0.4cm}   t,x \in \mathbb{R},
\end{align}
\begin{align}\label{1.4}
R_{2}^{\alpha,\beta}(t,x,\tilde{u},\tilde{v})=L_{2}^{\beta}(\tilde{v}(x,t))+\mathcal{M}(\tilde{u}(x,t), \tilde{v}(x,t))-g(x,t),  \hspace{0.4cm}   t,x \in \mathbb{R}
\end{align}
with initial conditions as:
\begin{align}\label{1.0a}
B_{1}\left( \tilde{u},\frac{\partial \tilde{u}}{\partial n}\right) =0, \hspace{0.4cm}  B_{2}\left( \tilde{v},\frac{\partial \tilde{v}}{\partial n}\right) = 0.
\end{align}
where $ \tilde{u} $, $\tilde{v}$ is the approximate solution to set of equations (\ref{1.1}) and (\ref{1.2}).
\end{dfn}
\begin{dfn} We define the sequences of functions $ \lbrace Q_{r}^{\alpha}(x,t)\rbrace_{r\in \mathbb{N}} $, $ \lbrace W_{h}^{\beta}(x,t)\rbrace_{h\in \mathbb{N}} $ as HP-sequences for equations (\ref{1.1}), (\ref{1.2}) respectively, where  $ Q_{r}^{\alpha}(x,t) $= $ \sum_{i=0}^{r}K_{r}^{i}\phi_{i} $, $ W_{h}^{\beta}(x,t)=\sum_{j=0}^{{h}}D_{h}^{j}\varphi_{j} $, where $ r,h\in \mathbb{N} $, $ K_{r}^{i}, D_{h}^{j}\in \mathbb{R}$. The functions involved in these sequences are called HP functions for (\ref{1.1}), (\ref{1.2}).
\end{dfn}
\begin{rem} If $$ \lim_{(r,h)\to (\infty, \infty)} R_{1}^{\alpha,\beta}(x,t,Q_{r}^{\alpha}(x,t),W_{h}^{\beta}(x,t))=0 \quad \text{and} \quad \lim_{(r,h)\to (\infty, \infty)} R_{2}^{\alpha,\beta}(x,t,Q_{r}^{\alpha}(x,t),W_{h}^{\beta}(x,t))=0,$$ then we say that HP-sequences $\lbrace Q_{r}^{\alpha}\rbrace $ and $\lbrace W_{h}^{\beta}\rbrace $ converge to the actual solution of coupled equations (\ref{1.1}), (\ref{1.2}).
\end{rem}
\begin{dfn} We call $ \tilde{u} $, $ \tilde{v} $ as the $ \epsilon $- approximate HP-solution to coupled equations (\ref{1.1}), (\ref{1.2}) on domain $\Omega$ if $$\mid R_{1}^{\alpha,\beta}(x,t,\tilde{u},\tilde{v})\mid<\epsilon \quad \text{and} \quad \mid R_{2}^{\alpha,\beta}(x,t,\tilde{u}, \tilde{v})\mid<\epsilon,$$
provided (\ref{1.0a}) is also satisfied by $ \tilde{u} $ and $\tilde{v}$.
\end{dfn}
\begin{dfn} We say that $ \tilde{u} $,  $ \tilde{v} $ constitute the weak $ \epsilon $-approximate HP-solution for set of equations (\ref{1.1}), (\ref{1.2}) on domain $\Omega $ if  $$\iint\limits_{\Omega} \ R^{2(\alpha,\beta)}_{1}(x,t,\tilde{u},\tilde{v})dx dt \leqslant \epsilon \quad \text{and} \quad  \iint\limits_{\Omega} \ R^{2(\alpha,\beta)}_{2}(x,t,\tilde{u},\tilde{v})dx dt \leqslant \epsilon$$
along with equation (\ref{1.0a}).
\end{dfn}
To calculate the constants $ K_{r}^{i} $ and $ D_{h}^{j} $ accurately, we propose the following steps:\\
\textbf{Step-I} We substitute the expression of $ \tilde{u} $, $ \tilde{v} $ in (\ref{1.3}) and (\ref{1.4}), and assume the transformed residual as:
\begin{align}\label{1.14}
\Re_{1}^{\alpha,\beta}(x,t,K_{r}^{i},D_{h}^{j})=R_{1}^{\alpha,\beta}(x,t,\tilde{u},\tilde{v})
\end{align}
and
\begin{align}\label{1.15}
\Re_{2}^{\alpha,\beta}(x,t,K_{r}^{i},D_{h}^{j})=R_{2}^{\alpha,\beta}(x,t,\tilde{u},\tilde{v}).
\end{align}
\textbf{Step-II} We consider the following functional:
\begin{equation}\label{1.16}
J(K_{r}^{i},D_{h}^{j})=\iint\limits_{\Omega} \left( \Re^{2(\alpha,\beta)}_{1}({x,t,K_{r}^{i},D_{h}^{j}})+ \Re^{2(\alpha,\beta)}_{2}({x,t,K_{r}^{i},D_{h}^{j}})\right) dx dt.
\end{equation}
Here, we calculate some constants $ K_{r}^{i} $ $ (K^{0}_{r},k^{1}_{r},\dots,K^{m}_{r}) $ and $ D_{h}^{j} $ $ (D^{0}_{h},D^{1}_{h},\dots,D^{n}_{h}) $ for $ m,n\in \mathbb{N} $ where $ m\leqslant r $, $ n\leqslant h $ as functions of remaining constants $ K_{r}^{i} $ ($ K_{r}^{m+1}, K_{r}^{m+2}, \dots,K_{r}^{r}$), $ D_{h}^{j} $ ($ D_{h}^{n+1}, D_{h}^{n+2}, \dots,D_{h}^{h}$) through given conditions (initial and boundary).
The remaining constants are determined through the minimization of functional (\ref{1.16}) (least square procedure).
\begin{rem} The constants can also be determined without using the initial /boundary conditions by directly minimizing (\ref{1.16}). If we represent these steps provide us the values of constants $ K_{r}^{i} $ as $ \tilde{K_{r}^{i}} $, and $ D_{h}^{j} $ as $ \tilde{D}_{h}^{j} $ then we can write HP-sequences as: $ Q_{r}^{\alpha}(x,t) $= $ \sum_{i=0}^{r}\tilde{K_{r}^{i}}\phi_{i} $ and $ W_{h}^{\beta}(x,t) $= $ \sum_{j=0}^{h}\tilde{D}_{h}^{j}\varphi_{j} $.
\end{rem}
\begin{lem} The property 
\begin{align*}
\lim_{(r,h)\to (\infty, \infty)}\iint\limits_{\Omega} \left( R^{2(\alpha,\beta)}_{1}(x,t,Q_{r}^{\alpha}(x,t),W_{h}^{\beta}(x,t))+R^{2(\alpha,\beta)}_{2}(x,t,Q_{r}^{\alpha}(x,t),W_{h}^{\beta}(x,t))\right) dx dt=0.
\end{align*}
is satisfied by  the HP-sequences $\lbrace{Q_{r}^{\alpha}(x,t)}\rbrace $ and $\lbrace{ W_{h}^{\beta}(x,t)}\rbrace  $ simultaneously. 
\end{lem}
\begin{proof}
Let $ Q_{r}^{\alpha}(x,t) $ and $ W_{h}^{\beta}(x,t) $ be two HP-functions developed on the basis of previous pertinent definitions.  Then following inequality is true:
\begin{align}\label{4.1}
\iint\limits_{\Omega} \left( R^{2(\alpha,\beta)}_{1}(x,t,Q_{r}^{\alpha}(x,t),W_{h}^{\beta}(x,t))+R^{2(\alpha,\beta)}_{2}(x,t,Q_{r}^{\alpha}(x,t),W_{h}^{\beta}(x,t))\right) dx dt\geqslant0.
\end{align}
Also from previous definitions, we have
\begin{align}
&\iint\limits_{\Omega} \left( R^{2(\alpha,\beta)}_{1}(x,t,u_{r}(x,t),v_{h}(x,t))+R^{2(\alpha,\beta)}_{2}(x,t,u_{r}(x,t), v_{h}(x,t))\right) dx dt& \nonumber \\ 
&\geqslant \iint\limits_{\Omega} \left( R^{2(\alpha,\beta)}_{1}(x,t,Q_{r}^{\alpha}(x,t),W_{h}^{\beta}(x,t))+R^{2(\alpha,\beta)}_{2}(x,t,Q_{r}^{\alpha}(x,t),W_{h}^{\beta}(x,t))\right) dx dt
\hspace{0.2cm} \forall r,h \in \mathbb{N}.&\label{4.2}
\end{align}
Combining inequalities (\ref{4.1}) and (\ref{4.2}), and taking $ \lim_{(r,h)\to (\infty, \infty)}$, we have 
\begin{align*}
&0 \leqslant \lim_{(r,h)\to (\infty, \infty)} \iint\limits_{\Omega} \left( R^{2(\alpha,\beta)}_{1}(x,t,Q_{r}^{\alpha}(x,t),W_{h}^{\beta}(x,t))+R^{2(\alpha,\beta)}_{2}(x,t,Q_{r}^{\alpha}(x,t),W_{h}^{\beta}(x,t))\right) dx dt\\
&\leqslant
\lim_{(r,h)\to (\infty, \infty)} \iint\limits_{\Omega} \left( R^{2(\alpha,\beta)}_{1}(x,t,u_{r}(x,t),v_{h}(x,t))
+R^{2(\alpha,\beta)}_{2}(x,t,u_{r}(x,t),v_{h}(x,t))\right) dx dt = 0.&
\end{align*}
Hence, we obtain the desired result.
\end{proof}
\begin{cor} The $ \epsilon $-approximate HP-solutions $ Q_{r}^{\alpha}(x,t) $ and $ W_{h}^{\beta}(x,t) $ are also the weak solutions to coupled equations (\ref{1.1}) and (\ref{1.2}). 
\end{cor}
\section{Numerical discussions of fractional non-linear wave equations}
In this section, we present three examples to illustrate the applicability of generalized least square homotopy perturbations (GLSHP). \\
\textbf{{Example 1.}}
One-dimensional space-fractional nonlinear wave equation is considered as \cite{Ghazanfari2011}:
\begin{align}\label{1.19}
\frac{\partial^{2 \alpha}u}{\partial x^{2\alpha}}-u\frac{\partial^{2}u}{\partial t^{2}} = 1-\frac{x^{2}+t^{2}}{2}\hspace{0.3cm} 0\leqslant x ,  t\leqslant 1, \hspace{0.2cm} 0<\alpha\leqslant 1,
\end{align}
with
\begin{align}\label{1.20}
u(t,0)= \dfrac{t^{2}}{2} , \hspace{0.4cm} \frac{\partial}{\partial x}u(t,0)= 0.
\end{align}
Exact solution to this example for $ \alpha=1 $ is $u(x,t)= \dfrac{x^{2}+t^{2}}{2}$ ( \cite{Ghazanfari2011}).\\
To solve equation (\ref{1.19}) with conditions (\ref{1.20}), we constitute homotopy equation as:
\begin{align}\label{1.22}
(1-p)\left[\frac{\partial^{2\alpha}U}{\partial x^{2\alpha}}- 1+\frac{x^{2}+t^{2}}{2}\right] +p\left[\frac{\partial^{2\alpha}U}{\partial x^{2\alpha}}-U\frac{\partial^{2} U}{\partial t^{2}}-1+ \frac{x^{2}+t^{2}}{2}\right] =0.
\end{align}
When $ p=0 $, equation (\ref{1.22}) gives $L^{2\alpha}[U(x,t,0)] = u_{0}(x,t),$ where initial guess $ u_{0}(x,t) $ will be determined from $\dfrac{\partial^{2\alpha}U}{\partial x^{2\alpha}}[u_{0}(x,t)]-1+\dfrac{x^{2}+t^{2}}{2} = 0, \hspace{0.3cm} u_{0}(t,0)=\dfrac{t^{2}}{2}, \hspace{0.3cm} \dfrac{\partial u_{0}}{\partial x}(t,0)=0$, whose solution is  $u_{0}=u_{0}^{(0)}(t,0)+u_{0}^{(1)}(t,0)x+J^{2\alpha}\left(\dfrac{2-x^{2}-t^{2}}{2}\right).$\\
Let expansion of $U(x,t,p)$ under homotopy considerations be
 \begin{align}\label{1.24}
U(x,t,p)=u_{0}(x,t)+pu_{1}(x,t)+p^{2}u_{2}(x,t)+p^{3}u_{3}(x,t)+\dots
\end{align}

\begin{figure}[htbp]
\subfloat[At  $\alpha = 0.89$]{\includegraphics[height=2.5in]{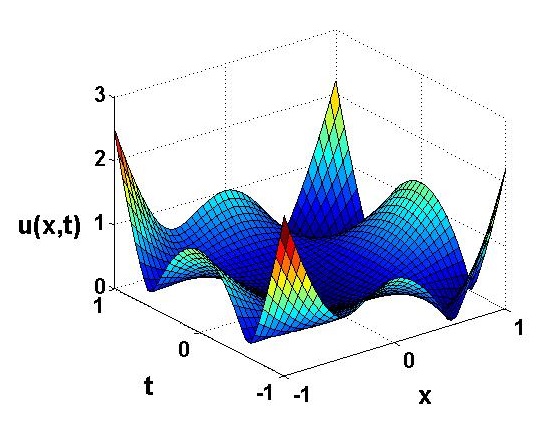}}
\subfloat[At $\alpha = 0.93$]{\includegraphics[height=2.5in]{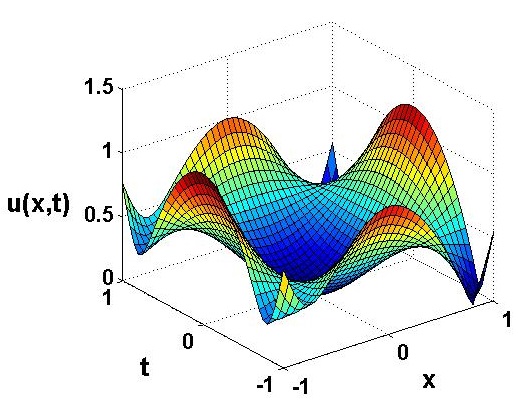}}\\
\subfloat[At $\alpha = 0.97$]{\includegraphics[height=2.5in]{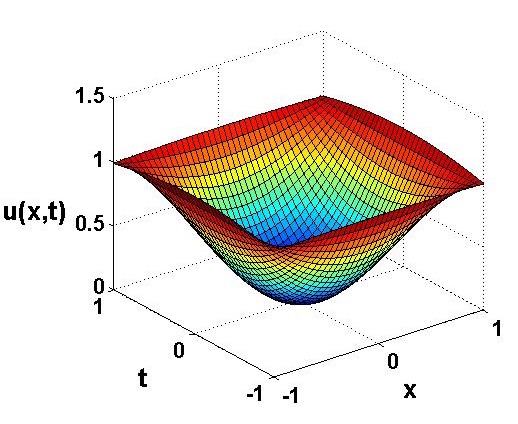}}
\subfloat[At $\alpha = 0.98$]{\includegraphics[height=2.5in]{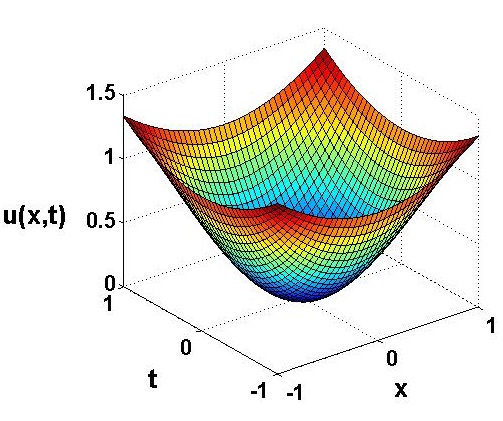}}\\
\subfloat[At $\alpha = 99$]{\includegraphics[height=2.4in]{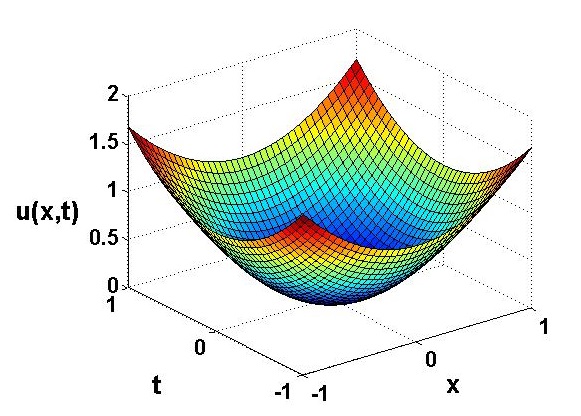}}
\subfloat[At $\alpha = 1$]{\includegraphics[height=2.4in]{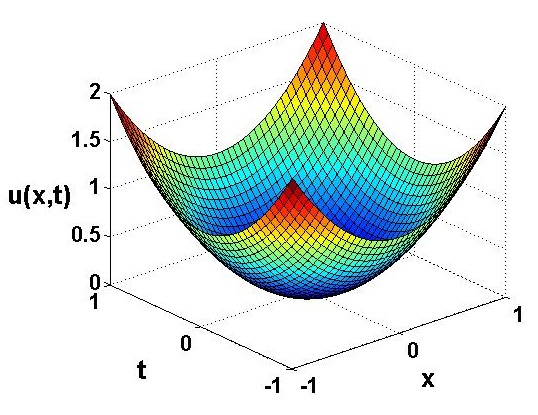}}
\caption{Plot of $\bar {u} (x,t)$ with respect to $ x, t $ at different values of $\alpha$.}\label{F3}
\end{figure}

\begin{figure}[htbp]
\subfloat[At $\beta = 0.89$]{\includegraphics[height=2.5in]{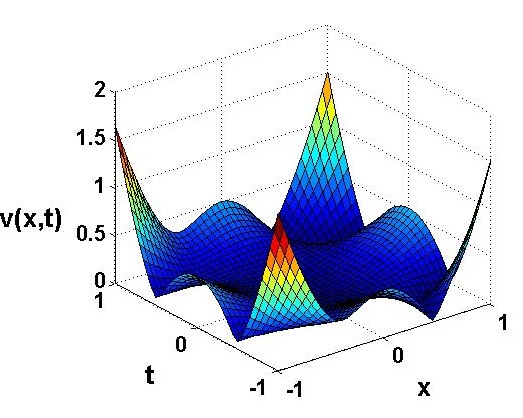}}
\subfloat[At $\beta = 0.93$]{\includegraphics[height=2.5in]{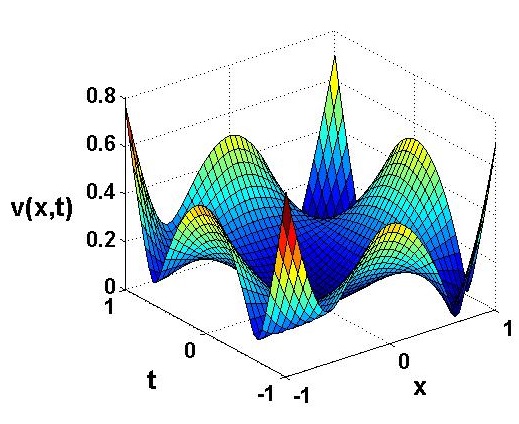}}\\
\subfloat[At $\beta = 0.97$]{\includegraphics[height=2.5in]{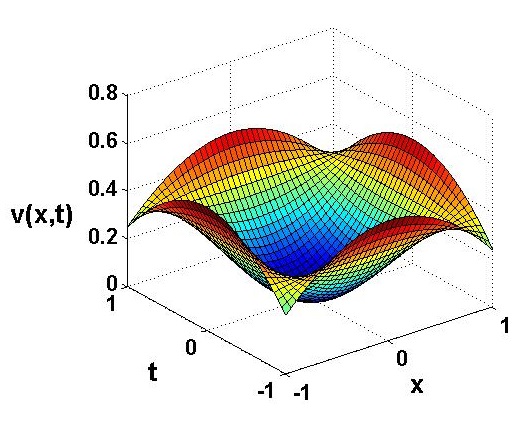}}
\subfloat[At $\beta = 0.98$]{\includegraphics[height=2.4in]{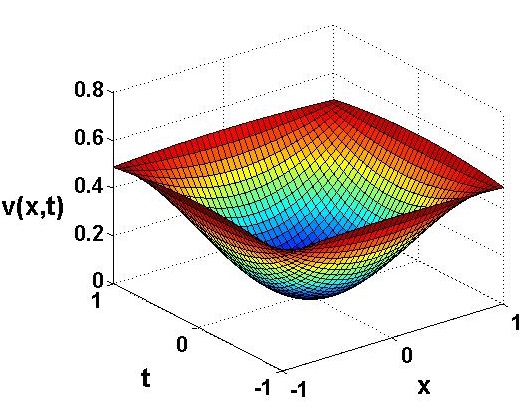}}\\
\subfloat[At $\beta = 0.99$]{\includegraphics[height=2.4in]{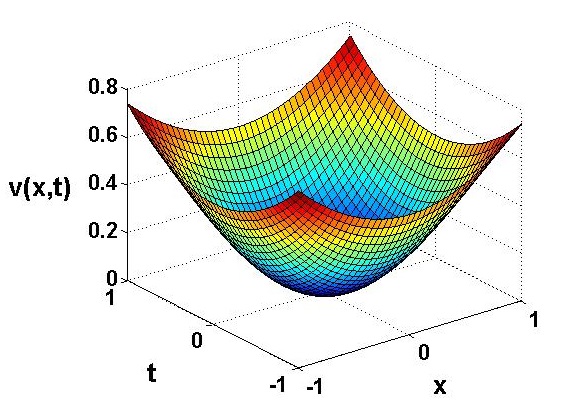}}
\subfloat[At $\beta = 1$]{\includegraphics[height=2.4in]{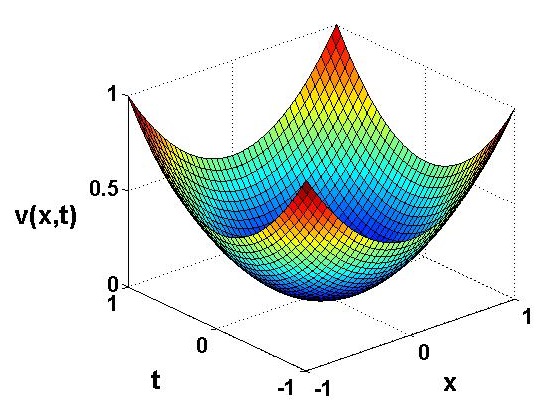}}
\caption{Plots of $\tilde{v} (x,t)$ with respect to $ x $ and $ t $ at different values of $ \beta $.}\label{F4}
\end{figure}

Substituting (\ref{1.24}) into (\ref{1.22}), we yield zero-$\text{th} $ order approximation as:
\begin{align*}
u_{0}(x,t)= \dfrac{t^{2}}{2}+\dfrac{x^{2\alpha}}{\Gamma(2\alpha+1)}-\dfrac{x^{2+2\alpha}}{\Gamma(2\alpha+3)}-\dfrac{t^{2}x^{2\alpha}}{2\Gamma(2\alpha+1)}
\end{align*}
Let $ Q_{0}=\left\lbrace \phi_1, \phi_2, \phi_3 \right\rbrace$ where $  \phi_1=\dfrac{x^{2\alpha}}{\Gamma(2\alpha+1)}$, $\phi_2=t^{2}$ and $\phi_3=\dfrac{x^{2\alpha}t^{2}}{\Gamma(2\alpha+1)}+\dfrac{x^{2+2\alpha}}{\Gamma(2\alpha+3)}$.\\
Here,
\begin{align*}
W^{\alpha}\left[ \phi_1, \phi_2, \phi_3 \right]  =
\begin{vmatrix}
\dfrac{x^{2\alpha}}{\Gamma(2\alpha+1)}& t^{2} &\left( \dfrac{x^{2\alpha}t^{2}}{\Gamma(2\alpha+1)} + \dfrac{x^{2+2\alpha}}{\Gamma(2\alpha+3)}\right)  \\ 
D^{\alpha} \dfrac{x^{2\alpha}}{\Gamma(2\alpha+1)} & D^{\alpha} t^{2} & D^{\alpha}\left( \dfrac{x^{2\alpha}t^{2}}{\Gamma(2\alpha+1)}+   \dfrac{x^{2+2\alpha}}{\Gamma(2\alpha+3)}\right)  \\
D^{2\alpha} \dfrac{x^{2\alpha}}{\Gamma(2\alpha+1)} & D^{2\alpha} t^{2} & D^{2\alpha} \left( \dfrac{x^{2\alpha}t^{2}}{\Gamma(2\alpha+1)} + \dfrac{x^{2+2\alpha}}{\Gamma(2\alpha+3)}\right)  \\
\end{vmatrix}
\end{align*}
For $\alpha=1$, $W^{1}\left[ \phi_1, \phi_2, \phi_3 \right]=-0.102\ne0$ when $x=0.2$ and $t=0.5$. Thus, by theorem on linear independence, functions $\phi_1$, $\phi_2$ and $\phi_3$ are L.I.\\
Now, let the zero-th order approximation to actual solution of (\ref{1.19}) and (\ref{1.20}) be
\begin{align}\label{1.32aa}
\tilde{u}=K_{0}\frac{x^{2\alpha}}{\Gamma(2\alpha+1)}+K_{1}t^{2}+K_{2}\left( t^{2}\dfrac{x^{2\alpha}}{\Gamma(2\alpha+1)}+\dfrac{x^{2+2\alpha}}{\Gamma(2\alpha+3)}\right) 
\end{align}
Let the residual function be 
\begin{align}\label{1.32aaa}
R^{\alpha}(x, t, \tilde{u})= \frac{\partial^{2\alpha}\tilde{u}}{\partial x^{2\alpha}}-\tilde{u}\frac{\partial^{2}\tilde{u}}{\partial t^{2}}+\frac{x^{2}+t^{2}-2}{2} 
\end{align}
and initial conditions are
\begin{align}\label{1.32aaa}
\tilde{u}(t,0)= \frac{t^{2}}{2} , \hspace{0.4cm} \frac{\partial}{\partial x}\tilde{u}(t,0)= 0
\end{align}
Using conditions (\ref{1.32aaa}) in equation (\ref{1.32aa}), we get $K_{1}=\dfrac{1}{2}$. Thus, the reduced approximate solution is
\begin{align}\label{1.34}
\tilde{u}=K_{0}\frac{x^{2\alpha}}{\Gamma(2\alpha+1)}+\frac{1}{2}t^{2}+K_{2}\left( t^{2}\dfrac{x^{2\alpha}}{\Gamma(2\alpha+1)}+\dfrac{x^{2+2\alpha}}{\Gamma(2\alpha+3)}\right) 
\end{align}
Modified residual ($ \Re^{\alpha}$) for this problem  can be written as\\
$ \Re^{\alpha}(x,t,K_{0},K_{2})= -1 + K_{0} + \dfrac{t^{2}}{2} + 
 K_{2} t^2 + \dfrac{x^{2}}{2}+ K_{2}\dfrac{x^{2}}{2} -\left( 1+ \dfrac{2 K_{2} x^{2\alpha}}{\Gamma(1+2\alpha)} \right)  \left( \dfrac{t^{2}}{2}+ 
 \dfrac{K_{0} x^{2\alpha}}{\Gamma(1+2\alpha)}+\left( \dfrac{K_{2} t^{2} x^{2\alpha}}{\Gamma(1 + 2\alpha)}+\dfrac{K_{2} x^{2 + 2\alpha}}{\Gamma(3 +2 \alpha)}\right)\right) $.\\
The functional $ J $ will be
\begin{equation*}
J(K_{0},K_{2})=\int_{0}^{1}\int_{0}^{1}\Re^{2(\alpha,\beta)}(x,t,K_{0},K_{2})dx dt.
\end{equation*}
Employing least square approximations (minimizing $J$), we receive two algebraic equations as: 
\begin{align*}
&\dfrac{1}{15 \Gamma(1 + 2 \alpha)^{4} \Gamma(3 + 2 \alpha)}\left\lbrace \left( \dfrac{40 K_{2}^{2} (3 K_{0}+K_{2}) \Gamma(3 + 2 \alpha)}{1 + 8 \alpha} \right) +4 
K_{2} \Gamma(1 + 2 \alpha)\left( \dfrac{30 K_{2}^{2}}{3 + 8 \alpha}+\dfrac{15 (2 K_{0} + K_{2}) \Gamma(3 + 2 \alpha)}{1 + 6 \alpha}\right) \right\rbrace &\\&
\dfrac{1}{15 \Gamma(1 + 2  \alpha)^{4} \Gamma(3 + 2 \alpha)}\left\lbrace \dfrac{5 \Gamma(1 + 
   2 \alpha)^{4} (-6 K_{2} + (3 + 2 \alpha) (-5 + 6 K_{0} + 3 K_{2}) \Gamma(3 + 2 \alpha)}{3 + 2 \alpha} \right\rbrace +\dfrac{1}{15 \Gamma(1 + 2  \alpha)^{2} \Gamma(3 + 2 \alpha)}&\\&
\left\lbrace \dfrac{40 K_{2}^{2}}{1 + 2 \alpha}+\dfrac{(10 (21 + 4 \alpha (5 - 7 K_{2}) - 15 K_{2}) K_{2} - 
   30 (3 + 4 \alpha) K_{0} (-1 + 4 K_{2})) \Gamma(3 + 2 \alpha)}{(1 + 4\alpha ) (3 + 4 \alpha)}  \right\rbrace -\dfrac{1}{15 \Gamma(1 + 2 \alpha) \Gamma(3 + 2 \alpha)}&\\
&\left\lbrace \dfrac{30 K_{2} (-1 + 2 K_{2})}{3 + 4 \alpha}+\dfrac{15 (-5 + 12 K_{0} + 7K_{2} +\alpha (-2 + 8 K_{0} + 6 K_{2})) \Gamma(3 + 2 \alpha)}{3 + 8 \alpha + 4 \alpha^{2}} \right\rbrace =0,
\end{align*}
\begin{align*}
&\dfrac{1}{60 \Gamma(1 + 2 \alpha)^{4} }\left( \dfrac{32 K_{2}^{2} (5 K_{0} + 3 K_{2}) }{1 + 8 \alpha}+\dfrac{32 K_{2} (15 K_{0}^{2} + 10 K_{0}K_{2}+ 3 K_{2}^{2}) }{1 + 8 \alpha}\right)+
\dfrac{1}{60 \Gamma(1 + 2 \alpha)^{4} \Gamma(3 + 2 \alpha)^{2}}&\\&
\left( \dfrac{48 K_{2} \Gamma(1 + 2 \alpha) \Gamma(
  3 + 2 \alpha) (10 (1 + 6 \alpha) K_{2} (2K_{0} + K_{2}) + (3 + 8 \alpha) (5 K_{0} + 
      4 K_{2}) )}{(1 + 6 \alpha) (3 + 8 \alpha)}\right)+\dfrac{16}{60 \Gamma(1 + 2 \alpha)^{3} \Gamma(3 + 2 \alpha)}&\\&
\left(\dfrac{10 K_{2}^{2} (3K_{0} + K_{2})}{3 + 8 \alpha}+\dfrac{3 (5 K_{0}^{2} + 5 K_{0} K_{2} + 2 K_{2}^{2}) \Gamma(3 + 2 \alpha)}{1 + 6 \alpha} \right)-\dfrac{1}{60 \Gamma(1 + 2 \alpha)^{4} \Gamma(3 + 2 \alpha)^{3}(3 + 4 \alpha) (5 + 4 \alpha) (3 + 8 \alpha + 4 \alpha^2)}&\\& 
\left\lbrace 4 \Gamma(1 + 2 \alpha)^{3} \Gamma(
  3 + 2 \alpha)^{2} (10 (3 + 8 \alpha + 4 \alpha^{2}) (3 (5 + 4 \alpha) K_{0} (-1 + 4 K_{2}) + 
    K_{2} (-62 + 57 K_{2} + 20 \alpha (-2 + 3 K_{2}))) \right\rbrace+&\\&
\dfrac{1}{60 \Gamma(1 + 2 \alpha)^{4} \Gamma(3 + 2 \alpha)^{3}}\left\lbrace (15 + 32 \alpha + 16 \alpha^{2}) (-50 + 105 K_{0} + 92 K_{2} + 
 \alpha (-20 + 90 K_{0} + 88 K_{2})) \Gamma(3 + 2 \alpha) \right\rbrace &\\& 
\dfrac{1}{15 \Gamma(1 + 2 \alpha)^{2} \Gamma(3 + 2 \alpha)^{2}}\left\lbrace \dfrac{240 K_{2}^{3}}{5 + 8 \alpha}+\dfrac{20 K_{2}^{2} \Gamma(3 + 2 \alpha)}{1 + 2 \alpha}+\dfrac{40 K_{2} (2 K_{0} +K_{2}) \Gamma(3 + 2 \alpha)}{1 + 2 \alpha} \right\rbrace-\dfrac{1}{15 \Gamma(1 + 2 \alpha)^{2} \Gamma(3 + 2 \alpha)^{2}}&\\&
\left\lbrace \dfrac{2 (30 (3 + 4 \alpha) K_{0}^{2} +K_{2} (-86 + 69 K_{2} + 44 \alpha (-2 + 3 K_{2})) + 
   5K_{0} (-21 + 30 K_{2}+ 4 \alpha (-5 + 14 K_{2}))) \Gamma(3 + 2 \alpha)^{2}}{(1 + 4 \alpha) (3 + 4 \alpha)} \right\rbrace +&\\&
\dfrac{1}{60\Gamma(3 + 2 \alpha)^{3}}\left\lbrace \dfrac{1440 (1 + 3 \alpha + 2 \alpha^{2}) K_{2}^{2}}{5 + 6 \alpha}+\dfrac{120 K_{2} \Gamma(3 + 2 \alpha)}{5 + 4 \alpha}-\dfrac{20 (-21 + 30 K_{0} + 38K_{2} + 2 \alpha (-3 + 6 K_{0} + 10 K_{2})) \Gamma(3 + 2 \alpha)^{2}}{(3 + 2 \alpha) (5 + 2 \alpha)}\right\rbrace&\\&
+\dfrac{1}{90} \left\lbrace (-71 + 90K_{0} + 65 K_{2})\right\rbrace=0. 
\end{align*}
Solving these equations simultaneously, we will get values of $K_0$ and $K_2$, and equation (\ref{1.32aa}) will be the required approximate solution.\\
In particular, when $ \alpha=1 $, we achieve $ K_{0}=1 $ and $K_{2}=0$, and thus $ \tilde{u}=\dfrac{x^{2}+t^{2}}{2}$ which is the exact solution \cite{Ghazanfari2011}.\\
\textbf{{Example 2.}} Consider the following space-fractional nonlinear wave equation \cite{Kaya1999}: 
\begin{align}\label{n1}
\frac{\partial^{2\alpha}u}{\partial x^{2\alpha}}-u\frac{\partial^{2}u}{\partial t^{2}}-2+2x^{2}+2t^{2}=0, \hspace{0.2cm} x,t\in [0,1],\hspace{0.2cm} 0<\alpha\leqslant 1
\end{align}
with initial conditions 
\begin{align}\label{n2}
u(t,0) = t^{2} ,\hspace{0.4cm}  \frac{\partial}{\partial x}u(t,0)=0,\hspace{0.4cm}u(0,x)=x^{2}.
\end{align}
For equation (\ref{n1}), homotopy perturbation equation (based on $ p $) is written as:
\begin{align}\label{n3}
(1-p)\left[ \frac{\partial^{2\alpha}U}{\partial x^{2\alpha}}-2+2x^{2}+2t^{2}\right] +p\left[\frac{\partial^{2\alpha}U}{2\partial x^{2\alpha}}-U\frac{\partial^{2}U}{\partial t^{2}}-2+2x^{2}+2t^{2} \right]=0.
\end{align}
In homotopy theory, for $ p=0 $ we get the initial fractional differential equation as:
$
L^{^{2\alpha}}[U(x,t,0)] = u_{0}(x,t)   
$
where initial solution $ u_{0}(x,t) $ is achieved by solving
$
 \dfrac{\partial^{2\alpha}U}{\partial x^{2\alpha}}[u_{0}(t,x)]-2+2x^{2}+2t^{2} = 0, \hspace{0.3cm} u_{0}(t,0)= t^{2},\hspace{0.3cm} \dfrac{\partial u_{0}}{\partial x}(t,0)= 0.
$
The solution of this fractional differential set up is readily obtained as  
$
u_{0}=u_{0}^{(0)}(t,0)+u_{0}^{(1)}(t,0)x+J^{2\alpha}[2-2x^{2}-2t^{2}].
$\\
Let us again take the expansion of $U(x,t,p)$ in the increasing powers of $p$ as 
\begin{align}\label{n5}
U(x,t,p)=u_{0}(t,x)+pu_{1}(t,x)+p^{2}u_{2}(t,x)+p^{3}u_{3}(t,x)+\dots
\end{align}
From equations (\ref{n5}) and (\ref{n3}), we have the following initially approximated solution:
\begin{align*}
u_{0}=t^{2}+\frac{2x^{2\alpha}}{\Gamma(2\alpha+1)}-\frac{2t^{2}x^{2\alpha}}{\Gamma(2\alpha+1)}-\frac{4x^{2+2\alpha}}{\Gamma(2\alpha+3)}.
\end{align*}
Let us consider the following set $Q_0=\lbrace{\phi_1, \phi_2, \phi_3\rbrace}$ where $\phi_1=t^2$, $\phi_2=\dfrac{x^{2\alpha}}{\Gamma (2\alpha+1)}$ and $\phi_3=\dfrac{t^{2}x^{2\alpha}}{\Gamma (2\alpha+1)}+\dfrac{x^{2+2\alpha}}{\Gamma (2\alpha+3)}.$\\
Since,
\begin{align*}
W^{\alpha}\left[ \phi_1, \phi_2, \phi_3\right] & =
\begin{vmatrix}
t^{2} &\dfrac{x^{2\alpha}}{\Gamma(2\alpha+1)}&\left( \dfrac{x^{2\alpha}t^{2}}{\Gamma(2\alpha+1)} + \dfrac{x^{2+2\alpha}}{\Gamma(2\alpha+3)}\right)  \\ 
D^{\alpha} t^{2}& D^{\alpha} \dfrac{x^{2\alpha}}{\Gamma(2\alpha+1)}  & D^{\alpha}\left( \dfrac{x^{2\alpha}t^{2}}{\Gamma(2\alpha+1)}+   \dfrac{x^{2+2\alpha}}{\Gamma(2\alpha+3)}\right)  \\
 D^{2\alpha} t^{2}&D^{2\alpha} \dfrac{x^{2\alpha}}{\Gamma(2\alpha+1)} & D^{2\alpha} \left( \dfrac{x^{2\alpha}t^{2}}{\Gamma(2\alpha+1)} + \dfrac{x^{2+2\alpha}}{\Gamma(2\alpha+3)}\right)
\end{vmatrix}&\\
&=t^{2}\dfrac{x^{\alpha}}{\Gamma(\alpha+1)}\left\lbrace \dfrac{6x^{\alpha}}{\Gamma(\alpha+1)}+4t+\dfrac{x^{2-\alpha}}{\Gamma(3-\alpha)} \right\rbrace-\dfrac{x^{2\alpha}}{\Gamma(2\alpha+1)}\left\lbrace \dfrac{12tx^{\alpha}}{\Gamma(\alpha+1)} + 8t^{2}+\dfrac{2x^{2-\alpha}t}{\Gamma(3-\alpha)}\right\rbrace.&
\end{align*}
Here, $ D^{\alpha} $=$\left(\dfrac{\partial}{\partial t}+\dfrac{\partial^{\alpha}}{\partial x^{\alpha}} \right)$ and $ D^{2\alpha} $ = $\left(\dfrac{\partial^{2}}{\partial t^{2}}+\dfrac{2\partial^{2}}{\partial t \partial x^{\alpha}}+\dfrac{\partial^{2\alpha}}{\partial x^{2\alpha}} \right)$.\\
At $ \alpha=1 $, $W^{1}\left[ \phi_1, \phi_2, \phi_3\right]=0.0444\ne0$ when $(x,t)=(0.3, 0.4).$\\
Thus by previous theorem, functions $\phi_1, \phi_2$ and $ \phi_3$ are linearly independent.
Now, if we consider
\begin{align}\label{n7}
\tilde{u}=K_{0}t^{2}+K_{1}\dfrac{x^{2\alpha}}{\Gamma (2\alpha+1)} +K_{2}\left( t^{2}\dfrac{x^{2\alpha}}{\Gamma (2\alpha+1)}+\dfrac{x^{2+2\alpha}}{\Gamma (2\alpha+3)}\right)
\end{align}
as the approximate solution to the given problem then residual $R^{\alpha}  $ with respect to $ \tilde{u} $ can be defined as
\begin{align}\label{n8}
R^{\alpha}(x, t, \tilde{u})= \frac{\partial^{2\alpha}\tilde{u}}{\partial x^{2\alpha}}-\tilde{u}\frac{\partial^{2}\tilde{u}}{\partial t^{2}}-2+2x^{2}+2t^{2}=0
\end{align}
with initial conditions 
\begin{align}\label{n10}
\tilde{u}(t,0) = t^{2} ,\hspace{0.4cm}  \frac{\partial}{\partial x}\tilde{u}(t,0)=0.
\end{align}
Utilizing the conditions (\ref{n10}) in equation (\ref{n7}), we get $ K_{0}=1 $.\\
Now, the modified remainder $ \Re^{\alpha}$ becomes
\begin{align*}
&\Re_{1}^{\alpha}(x,t,K_{1},K_{2})=-2 + K_{1} + 2 t^{2} +K_{2} t^{2} + 2 x^{2} +\dfrac{ K_{2} x^{2}}{2}- \left(  2+\dfrac{ 2 K_{2} x^{2\alpha}}{\Gamma(1 + 2\alpha)} \right) \left(t^{2} + \dfrac{K_{1}x^{2\alpha}}{ \Gamma(1 + 2\alpha)} +\dfrac{K_{2} t^{2} x^{2\alpha}}{\Gamma(1 + 2\alpha)}+\dfrac{K_{2} x^{2 + 2\alpha}}{ \Gamma(3 + 2\alpha)} \right). 
\end{align*}
Further, to find the values of remaining parameters $K_1$ and $K_2$, we compose the functional $\textit{J}$ as  
\begin{align*}
J(K_{1},K_{2})=\int_{0}^{1}\int_{0}^{1}\Re^{2\alpha}(x,t,K_{1},K_{2})dx dt,
\end{align*}
and compute the values of parameters by minimizing this functional through least square procedure which gives two two equations of algebraic nature as:
\begin{align*}
&\dfrac{1}{15 \Gamma(1 + 2 \alpha)^{4} \Gamma(3 + 2 \alpha)}\left\lbrace \left( \dfrac{40 K_{2}^{2} (3 K_{1}+K_{2}) \Gamma(3 + 2 \alpha)}{1 + 8 \alpha} \right) +120K_{2}\Gamma(1 + 2 \alpha) 
\left( \dfrac{ K_{2}^{2}}{3 + 8 \alpha}+\dfrac{15 (2 K_{1} + K_{2}) \Gamma(3 + 2 \alpha)}{1 + 6 \alpha}\right) \right\rbrace &\\&
\dfrac{1}{15 \Gamma(1 + 2  \alpha)^{4} \Gamma(3 + 2 \alpha)}\left\lbrace \dfrac{5 \Gamma(1 + 
   2 \alpha)^{4} (-12 K_{2} + (3 + 2 \alpha) (-8 + 6 K_{1} + 3 K_{2}) \Gamma(3 + 2 \alpha)}{3 + 2 \alpha} \right\rbrace +\dfrac{1}{15 \Gamma(1 + 2  \alpha)^{2} \Gamma(3 + 2 \alpha)}&\\&
\left\lbrace \dfrac{40 K_{2}^{2}}{1 + 2 \alpha}+\dfrac{(-60 (3 + 4 \alpha) K_{1} (-1 + K_{2}) + 
   5 (48 + 4\alpha (8 - 7 K_{2}) - 15 K_{2}) K_{2}) \Gamma(3 + 2 \alpha)}{(1 + 4\alpha ) (3 + 4 \alpha)}  \right\rbrace -\dfrac{1}{15 \Gamma(1 + 2 \alpha) \Gamma(3 + 2 \alpha)}&\\&
   \left\lbrace \dfrac{30 K_{2} (-2 + 2 K_{2})}{3 + 4 \alpha}+\dfrac{15 (-8 + 4 (3 + 2 \alpha) K_{1} + (7 + 6 \alpha) K_{2})  \Gamma(3 + 2 \alpha)}{3 + 8 \alpha + 4 \alpha^{2}} \right\rbrace =0, 
\end{align*}
\begin{align*}
&\dfrac{1}{60 \Gamma(1 + 2 \alpha)^{4} }\left( \dfrac{32 K_{2}^{2} (5 K_{1} + 3 K_{2}) }{1 + 8 \alpha}+\dfrac{32 K_{2} (15 K_{1}^{2} + 10 K_{1}K_{2}+ 3 K_{2}^{2}) }{1 + 8 \alpha}\right)+
\dfrac{1}{60 \Gamma(1 + 2 \alpha)^{4} \Gamma(3 + 2 \alpha)^{2}}&\\&
\left( \dfrac{96 K_{2} \Gamma(1 + 2 \alpha) \Gamma(
  3 + 2 \alpha) (5 (1 + 6 \alpha) K_{2} (2K_{1} + K_{2}) + (3 + 8 \alpha) (5 K_{1} + 
      4 K_{2}) )}{(1 + 6 \alpha) (3 + 8 \alpha)}\right)+\left( \dfrac{1}{60 \Gamma(1 + 2 \alpha)^{4} \Gamma(3 + 2 \alpha)^{3}}\right) &\\&
 \left( \dfrac{1}{(3 + 4 \alpha) (5 + 4 \alpha) (3 + 8 \alpha + 4 \alpha^{2})} \right)     
8 \Gamma(1 + 2 \alpha)^{3} \Gamma(
  3 + 2 \alpha)^{2} (5 (3 + 8 \alpha + 4 \alpha^{2}) (12 (5 + 4 \alpha) K_{1} (-1 + K_{2}) + &\\&
 K_{2} (-128 + 57 K_{2} + \alpha (-64 + 60 K_{2}))) + (15 + 32 \alpha + 
      16 \alpha^{2}) (-80 + 15 (7 + 6 \alpha) K_{1}+ (92 + 88 \alpha) K_{2}) \Gamma(3 + 2 \alpha))&\\&      
\dfrac{32}{60 \Gamma(1 + 2 \alpha)^{3}\Gamma(3 + 2 \alpha) }\left(\dfrac{5 K_{2}^{2} (3K_{1} + K_{2})}{3 + 8 \alpha}+\dfrac{3 (5 K_{1}^{2} + 5 K_{1} K_{2} + 2 K_{2}^{2}) \Gamma(3 + 2 \alpha)}{1 + 6 \alpha} \right)-\dfrac{8}{60 \Gamma(1 + 2 \alpha)^{2} \Gamma(3 + 2 \alpha)^{2}}&\\&
\left( \dfrac{240 K_{2}^{3}}{5 + 8 \alpha}+\dfrac{20 K_{2}^{2} \Gamma(3 + 2 \alpha)}{1 + 2 \alpha} +\dfrac{40 K_{2} (2 K_{1} +K_{2}) \Gamma(3 + 2 \alpha)}{1 + 2 \alpha}\right)-\dfrac{1}{60 \Gamma(1 + 2 \alpha)^{4} \Gamma(3 + 2 \alpha)^{3}(1 + 4 \alpha) (3 + 4 \alpha)}  &\\&
\left\lbrace (30 (3 + 4 \alpha) K_{1}^{2}+10 K_{1}(3 (-8 + 5K_{2}) + 4 \alpha(-4 + 7 K_{2})) + 
  K_{2} (-224 + 69 K_{2} + 12 \alpha (-16 + 11 K_{2}))) \Gamma(3 + 2 \alpha)^{2} \right\rbrace&\\&
\dfrac{1}{60\Gamma(3 + 2 \alpha)^{3}}\left\lbrace \dfrac{2880 (1 + 3 \alpha + 2 \alpha^{2}) K_{2}^{2}}{5 + 6 \alpha}+\dfrac{480 K_{2} \Gamma(3 + 2 \alpha)}{5 + 4 \alpha}-\dfrac{80 (-12 + 3 (5 + 2 \alpha) K_{1}+ (19 + 10 \alpha) K_{2})) \Gamma(3 + 2 \alpha)^{2}}{(3 + 2 \alpha) (5 + 2 \alpha)}\right\rbrace&\\&
+\dfrac{1}{90} \left\lbrace (-104 + 90 K_{1} + 65 K_{2})\right\rbrace=0. 
\end{align*}
The solution of these two equations when combined with the value of  parameter $K_0 $ and equation $(\ref{n7})$ will produce the required solution.\\
Further when $ \alpha=1 $, we have  $\tilde{u}=t^{2}+x^{2}$, the exact solution \cite{Kaya1999}. \\
Now to portray the clear picture of the scheme, we take the following example which consists of coupled system of equations.\\
\textbf{{Example 3.}} Let us consider the following nonlinear space-fractional system of equations as \citep{Biazar2013}:
\begin{align}\label{2.14a}
\frac{\partial^{2\alpha}u}{\partial x^{2\alpha}}-v\frac{\partial^{2}u}{\partial t^{2}}-u\frac{\partial^{2}v}{\partial t^{2}}-2+2x^{2}+2t^{2}=0, \hspace{0.3cm} 0<\alpha\leqslant 1
\end{align}
\begin{align}\label{1.38}
\frac{\partial^{2\beta}v}{\partial x^{2\beta}}-v\frac{\partial^{2}v}{\partial t^{2}}+u\frac{\partial^{2}u}{\partial t^{2}}-1-\frac{3}{2}x^{2}-\frac{3}{2}t^{2}=0, \hspace{0.3cm}x,t\in [0,1],\hspace{0.2cm} 0< \beta\leqslant 1
\end{align} 
with
\begin{align}\label{1.39}
u(t,0) = t^{2} ,\hspace{0.4cm}  \frac{\partial}{\partial x}u(t,0)=0,
\end{align}
\begin{align}\label{1.40}
v(t,0) = \frac{t^{2}}{2} , \hspace{0.4cm}  \frac{\partial}{\partial x}v(t,0)=0. 
\end{align}
Writing the homotopy equations for equations (\ref{2.14a}) and (\ref{1.38}), we have:
\begin{align}\label{1.43}
(1-p)\left[ \frac{\partial^{2\alpha}U}{\partial x^{2\alpha}}-2+2x^{2}+2t^{2}\right] +p\left[\frac{\partial^{2\alpha}U}{\partial x^{2\alpha}}-V\frac{\partial^{2}U}{\partial t^{2}}-U\frac{\partial^{2}V}{\partial t^{2}}-2+2x^{2}+2t^{2} \right]=0,
\end{align}
\begin{align}\label{1.44}
(1-p)\left[ \frac{\partial^{2\beta}V}{\partial x^{2\beta}}-1-\frac{3}{2}x^{2}-\frac{3}{2}t^{2}\right] +p\left[\frac{\partial^{2\beta}V}{\partial x^{2\beta}}-V\frac{\partial^{2}V}{\partial t^{2}}-U\frac{\partial^{2}U}{\partial t^{2}}-1-\frac{3}{2}x^{2}-\frac{3}{2}t^{2} \right]=0. 
\end{align}
When $ p=0 $, these two equations [$(\ref{1.43})$, $(\ref{1.44})$] give $L^{^{2\alpha}}[U(x,t,0)] = u_{0}(x,t)  $ and $ L^{2\beta}[V(x,t,0)] = v_{0}(x,t) $, where initial approximations $ u_{0}(x,t) $ and $v_{0}(x,t)  $ can be determined from following respective equations:
\begin{align}\label{1.45aa}
 L^{2\alpha}[u_{0}(x,t)]-2+2x^{2}+2t^{2} = 0, \hspace{0.3cm} u_{0}(0,t)= t^{2},\hspace{0.3cm} \frac{\partial u_{0}}{\partial x}(0,t)= 0,
\end{align}
\begin{align}\label{1.46}
L^{2\beta}[v_{0}(x,t)]-1-\frac{3}{2}x^{2}-\frac{3}{2}t^{2}= 0, \hspace{0.3cm} v_{0}(0,t)= \frac{t^{2}}{2},\hspace{0.3cm} \frac{\partial v_{0}}{\partial x}(0,t)= 0.
\end{align}
The solution to equations (\ref{1.45aa}) and (\ref{1.46}) can be easily obtained from the basic theory of fractional calculus as
\begin{align*}
 u_{0}(x,t)=u_{0}^{(0)}(0,t)+x u_{0}^{(1)}(0,t)+J^{2\alpha}\lbrace2-2x^{2}-2t^{2}\rbrace,
\end{align*}
\begin{align*}
v_{0}(x,t)=v_{0}^{(0)}(0,t)+x v_{0}^{(1)}(0,t)+J^{2\beta}\lbrace1+\frac{3}{2}x^{2}+\frac{3}{2}t^{2}\rbrace.
\end{align*}
Following homotopy theory, and assuming the expansion of $U(x,t,p)  $ and $V(x,t,p)  $ in terms of increasing powers of $ p $ as 
\begin{align}\label{1.47}
U(x,t,p)=u_{0}(x,t)+pu_{1}(x,t)+p^{2}u_{2}(x,t)+p^{3}u_{3}(x,t)+...
\end{align}
\begin{align}\label{1.48}
V(x,t,p)=v_{0}(x,t)+pv_{1}(x,t)+p^{2}v_{2}(x,t)+p^{3}v_{3}(x,t)+...
\end{align}
It is a well-known phenomenon in homotopy (for equations (\ref{1.47}) and (\ref{1.48})) that $ U(x,t,p) \rightarrow u(x,t) $ and $V(x,t,p) \rightarrow v(x,t) $ when $ p$ varies from $0$ to 1.\\
Through equations (\ref{1.43}),(\ref{1.44}) and equations (\ref{1.47}), (\ref{1.48}), we achieve the initial solution as
\begin{align}\label{1.53}
 u_{0}=t^{2}+\frac{2x^{2\alpha}}{\Gamma(2\alpha+1)}-\frac{4x^{2+2\alpha}}{\Gamma(2\alpha+3)}-\frac{2t^{2}x^{2\alpha}}{\Gamma(2\alpha+1)},
\end{align}
\begin{align}\label{1.54}
 v_{0}=\frac{t^{2}}{2}+\frac{x^{2\beta}}{\Gamma(2\beta+1)}+\frac{3x^{2+2\beta}}{\Gamma(3+2\beta)}+\frac{3t^{2}x^{2\beta}}{2\Gamma(2\beta+1)}.
\end{align}
Following \textit{generalized least square homotopy perturbations}, we compose two function sets as:
$Q_{0} =\left\lbrace \phi_1, \phi_2, \phi_3 \right\rbrace $ and $S_{0} =\left\lbrace \varphi_1, \varphi_2, \varphi_3 \right\rbrace $ where $ \phi_1=\dfrac{x^{2\alpha}}{\Gamma(2\alpha+1)}$, $ \phi_2=t^{2}$, $ \phi_3=\dfrac{t^{2}x^{2\alpha}}{\Gamma(2\alpha+1)}+\dfrac{x^{2+2\alpha}}{\Gamma(2\alpha+3)}$, $ \varphi_1=\dfrac{x^{2\beta}}{\Gamma(2\beta+1)}$, $ \varphi_2=t^{2}$ and $ \varphi_3=\dfrac{t^{2}x^{2\beta}}{\Gamma(2\beta+1)}+\dfrac{x^{2+2\beta}}{\Gamma(3+2\beta)}.$\\
Here, we have two partial fractional Wronskians $W^{\alpha}_{1}$ and $W^{\beta}_{2}$ which are defined as:
\[
W^{\alpha}_{1}\left[\phi_1, \phi_2, \phi_3\right]  =
\begin{vmatrix}
 \dfrac{x^{2\alpha}}{\Gamma(2\alpha+1)}&t^{2}&\left( \dfrac{x^{2\alpha}t^{2}}{\Gamma(2\alpha+1)} + \dfrac{x^{2+2\alpha}}{\Gamma(2\alpha+3)}\right)  \\ 
 D^{\alpha} \dfrac{x^{2\alpha}}{\Gamma(2\alpha+1)}&D^{\alpha} t^{2}  & D^{\alpha}\left( \dfrac{x^{2\alpha}t^{2}}{\Gamma(2\alpha+1)}+   \dfrac{x^{2+2\alpha}}{\Gamma(2\alpha+3)}\right)  \\
D^{2\alpha} \dfrac{x^{2\alpha}}{\Gamma(2\alpha+1)} &D^{2\alpha} t^{2}& D^{2\alpha} \left( \dfrac{x^{2\alpha}t^{2}}{\Gamma(2\alpha+1)} + \dfrac{x^{2+2\alpha}}{\Gamma(2\alpha+3)}\right)  \\
\end{vmatrix}
\]
and 
\[
W^{\beta}_{2}\left[\varphi_1, \varphi_2, \varphi_3\right]  =
\begin{vmatrix}
 \dfrac{x^{2\beta}}{\Gamma(2\alpha+1)}&t^{2}&\left( \dfrac{x^{2\beta}t^{2}}{\Gamma(2\beta+1)} + \dfrac{x^{2+2\beta}}{\Gamma(2\beta+3)}\right)  \\ 
{D}^{\beta} \dfrac{x^{2\beta}}{\Gamma(2\beta+1)}&{D}^{\beta} t^{2}  & {D}^{\beta}\left( \dfrac{x^{2\beta}t^{2}}{\Gamma(2\beta+1)}+   \dfrac{x^{2+2\beta}}{\Gamma(2\beta+3)}\right)  \\
{D}^{2\beta} \dfrac{x^{2\beta}}{\Gamma(2\beta+1)} &{D}^{2\beta} t^{2}& {D}^{2\beta} \left( \dfrac{x^{2\beta}t^{2}}{\Gamma(2\beta+1)} + \dfrac{x^{2+2\beta}}{\Gamma(2\beta+3)}\right) \\
\end{vmatrix}.
\] 
Here $ D^{\alpha} $ = $\left(\dfrac{\partial}{\partial t}+\dfrac{\partial^{\alpha}}{\partial x^{\alpha}} \right)$, $ D^{2\alpha} $ = $\left(\dfrac{\partial^{2}}{\partial t^{2}}+\dfrac{2\partial^{2}}{\partial t \partial x^{\alpha}}+\dfrac{\partial^{2\alpha}}{\partial x^{2\alpha}} \right)$, $ {D}^{\beta}$ = $ \left(\dfrac{\partial}{\partial t}+\dfrac{\partial}{\partial x^{\beta}} \right) $ and  $ {D}^{2\beta} $ = $\left(\dfrac{\partial^{2}}{\partial t^{2}}+\dfrac{2\partial^{2}}{\partial t \partial x^{\beta}}+\dfrac{\partial^{2\beta}}{\partial x^{2\beta}} \right)$.\\
The simplification to these Wronskian gives
\begin{align*}
W_{1}^{\alpha}\left[\phi_1, \phi_2, \phi_3\right]  =\dfrac{x^{2\alpha}}{\Gamma(2\alpha+1)}\left\lbrace \dfrac{12tx^{\alpha}}{\Gamma(\alpha+1)} + 8t^{2}+\dfrac{2x^{2-\alpha}t}{\Gamma(3-\alpha)}\right\rbrace-t^{2}\dfrac{x^{\alpha}}{\Gamma(\alpha+1)}\left\lbrace \dfrac{6x^{\alpha}}{\Gamma(\alpha+1)}+4t+\dfrac{x^{2-\alpha}}{\Gamma(3-\alpha)} \right\rbrace  
\end{align*}
\begin{align*}
W_{2}^{\beta}\left[\varphi_1, \varphi_2, \varphi_3\right] =\dfrac{x^{2\beta}}{\Gamma(2\beta+1)}\left\lbrace \dfrac{12tx^{\beta}}{\Gamma(\beta+1)} + 8t^{2}+\dfrac{2x^{2-\beta}t}{\Gamma(3-\beta)}\right\rbrace-t^{2}\dfrac{x^{\beta}}{\Gamma(\beta+1)}\left\lbrace \dfrac{6x^{\beta}}{\Gamma(\beta+1)}+4t+\dfrac{x^{2-\beta}}{\Gamma(3-\beta)} \right\rbrace  
\end{align*}
For $ \alpha=1 $ and $ \beta=1 $, we have
$
  W_{1}^{1}\left[\phi_1, \phi_2, \phi_3\right]=7tx^{3}-3t^{2}x^{2}-4t^{3}x  
$ and $  W_{2}^{1}\left[\varphi_1, \varphi_2, \varphi_3\right]=7tx^{3}-3t^{2}x^{2}-4t^{3}x.  $\\
As we know that $ x,t\in [0,1] $, therefore, if we choose $( x,t)=(0.2,0.5) $ then
$
W_{1}^{1}\left[\phi_1, \phi_2, \phi_3\right]=-0.102\neq 0$ and $ W_{2}^{1} \left[\varphi_1, \varphi_2, \varphi_3\right]=-0.102\neq 0$. Thus the functions chosen in sets $Q_0$ and $S_0$ are linearly independent in vector space structure.\\
Based on these sets, the $ 0^{th} $- order approximate solution to present example is assumed as:
\begin{align}\label{1.57}
\tilde{u}=K_{0}\dfrac{x^{2\alpha}}{\Gamma(2\alpha+1)}+K_{1}t^{2}+K_{2}\left( \dfrac{t^{2}x^{2\alpha}}{\Gamma(2\alpha+1)}+\dfrac{x^{2+2\alpha}}{\Gamma(2\alpha+3)}\right), 
\end{align}
\begin{align}\label{1.58}
\tilde{v}=D_{0}\dfrac{x^{2\beta}}{\Gamma(2\beta+1)}+D_{1}{t^{2}}+D_{2}\left( \dfrac{t^{2}x^{2\beta}}{\Gamma(2\beta+1)}+\dfrac{x^{2+2\beta}}{\Gamma(3+2\beta)}\right). 
\end{align}
Now the residuals $R_{1}^{\alpha,\beta}  $ and $ R_{2}^{\alpha,\beta} $ depending on $ \tilde{u} $ and $\tilde{v}  $ are presented as:
\begin{align}\label{1.56a}
R_{1}^{\alpha,\beta}(x, t, \tilde{u},\tilde{v})= \frac{\partial^{2\alpha}\tilde{u}}{\partial x^{2\alpha}}-\tilde{v}\frac{\partial^{2}\tilde{u}}{\partial t^{2}}-\tilde{u}\frac{\partial^{2}\tilde{v}}{\partial t^{2}}-2+2x^{2}+2t^{2},
\end{align}
\begin{align}\label{1.56aa}
R_{2}^{\alpha,\beta}(x, t, \tilde{u},\tilde{v})=\frac{\partial^{2\beta}\tilde{v}}{\partial x^{2\beta}}-\tilde{v}\frac{\partial^{2}\tilde{v}}{\partial t^{2}}+\tilde{u}\frac{\partial^{2}\tilde{u}}{\partial t^{2}}-1-\frac{3}{2}x^{2}-\frac{3}{2}t^{2}
\end{align}
with
\begin{align}\label{1.56aaa}
\tilde{u}(t,0) = t^{2} ,\hspace{0.4cm}  \frac{\partial}{\partial x}\tilde{u}(t,0)=0,
\end{align}
\begin{align}\label{1.56aaaa}
\tilde{v}(t,0) = \frac{t^{2}}{2} , \hspace{0.4cm}  \frac{\partial}{\partial x}\tilde{v}(t,0)=0. 
\end{align}
Substituting the conditions (\ref{1.56aaa}) and (\ref{1.56aaaa}) into equations (\ref{1.57}) and (\ref{1.58}), we receive $ K_{1}=1 $ and $ D_{1}=\frac{1}{2} $ which further gives 
\begin{align}\label{1.59}
\tilde{u}=K_{0}\dfrac{x^{2\alpha}}{\Gamma(2\alpha+1)}+t^{2}+K_{2}\left( \dfrac{t^{2}x^{2\alpha}}{\Gamma(2\alpha+1)}+\dfrac{x^{2+2\alpha}}{\Gamma(2\alpha+3)}\right),
\end{align}
\begin{align}\label{1.60}
\tilde{v}=D_{0}\dfrac{x^{2\beta}}{\Gamma(2\beta+1)}+\dfrac{t^{2}}{2}+D_{2}\left( \dfrac{t^{2}x^{2\beta}}{\Gamma(2\beta+1)}+\dfrac{x^{2+2\beta}}{\Gamma(3+2\beta)}\right) .
\end{align}
Now, revised residuals $ \Re_{1}^{\alpha,\beta}$ and $ \Re_{2}^{\alpha,\beta}$ to the problem are
\begin{align*}
&\Re_{1}^{\alpha,\beta}(x,t,K_{0},K_{2},D_{0},D_{2})=-2 + K_{0} + 2 t^{2} + K_{2} t^{2} + 2 x^{2} + \dfrac{K_{2}x^{2}}{2} -  \left( t^{2} +\dfrac{K_{0} x^{2\alpha}}{\Gamma(1 + 2\alpha)} + \dfrac{K_{2} t^{2} x^{2\alpha}}{\Gamma(1 + 2\alpha)} + \dfrac{K_{2} x^{2 + 2\alpha}}{\Gamma(3 + 2\alpha)}\right) &\\&  
\left( 1+\dfrac{ 2 D_{2} x^{2\beta}}{\Gamma(1 +2 \beta)}\right)  -\left( 2+ \dfrac{2 K_{2}x^{2\alpha}}{\Gamma(1 + 2\alpha)}\right) \left( \dfrac{t^{2}}{2}+ D_{0} \dfrac{x^{2\beta}}{\Gamma(1 + 2\beta)} + D_{2} t^{2} \dfrac{x^{2\beta}}{\Gamma(1 + 2\beta)} + \dfrac{D_{2}x^{2 + 2\beta}}{
    \Gamma(3 + 2\beta)}\right), 
\end{align*}
\begin{align*}
&\Re_{2}^{\alpha,\beta}(x,t,K_{0},K_{2},D_{0},D_{2})=-1+D_{0} - \dfrac{3 t^{2}}{2} + D_{2} t^{2} - 
\dfrac{3 x^{2}}{2} +\dfrac{D_{2}x^{2}}{2} + \left( 2 + 2 K_{2}\dfrac{x^{2\alpha}}{\Gamma(1 + 2\alpha)} \right)  &\\&\left( t^{2} +\dfrac{K_{0} x^{2\alpha}}{ \Gamma(1 + 2\alpha)}  +\dfrac{K_{2} t^{2} x^{2\alpha}}{\Gamma(1 + 2\alpha)}  + \dfrac{K_{2} x^{2 + 2\alpha}}{ \Gamma(3 + 2\alpha)} \right)
  - \left( 1 + \dfrac{2 D_{2} x^{2\beta}}{\Gamma(1 + 2\beta)}\right)
   \left(  \dfrac{t^{2}}{2} + \dfrac{D_{0} x^{2\beta}}{ \Gamma(1 + 2\beta)} +\dfrac{D_{2} t^{2} x^{2\beta}}{\Gamma(1 + 2\beta)}  + \dfrac{D_{2}x^{2 + 2\beta}}{ \Gamma(3 + 2\beta)}\right) .
\end{align*}
The values of remaining parameters are now determined through the functional $\textit{J}$ 
\begin{align*}
J(K_{0},K_{2},D_{0},D_{2})=\int_{0}^{1}\int_{0}^{1}(\Re^{2(\alpha,\beta)}_{1}(x,t,K_{0},K_{2},D_{0},D_{2})+\Re^{2(\alpha,\beta)}_{2}(x,t,K_{0},K_{2},D_{0},D_{2}))dx dt
\end{align*}
and its minimization in the form of algebraic equations 
\begin{align*}
\left\lbrace \dfrac{\partial J}{\partial K_{0}}=0, \dfrac{\partial J}{\partial K_{2}}=0, \dfrac{\partial J}{\partial D_{0}}=0, \dfrac{\partial J}{\partial D_{2}}=0 \right\rbrace.
\end{align*}
Here, detailed form of equations has been ignored to conserve space. Thus, required solution will be obtained by putting these parametric values in equations (\ref{1.57}) and (\ref{1.58}).\\
Moreover, we get values of constants as $ K_{0}=2 $, $K_{2}=0$ , $ D_{0}=1 $ and $ D_{2}= 0$ for $ (\alpha, \beta)=(1,1) $ which further gives 
$\tilde{u}=x^{2}+t^{2}$ and $\tilde{v}=0.5x^{2}+\dfrac{t^{2}}{2}$ as exact solution \citep{Biazar2013}.
Figures (\ref{F3}) and (\ref{F4}) exhibit the dynamics of  $\tilde{u}(x,t)$, $\tilde{v}(x,t)$ at different values of $ \alpha $ and $ \beta $ with respect to various $ x $ and $ t $.

\section{Conclusions}
In this paper, an effective combination of fractional homotopy perturbations and least square approximations have been constructed to obtain reliable solutions for fractional partial differential equations. A novel idea of fractional partial Wronskian is also developed to verify the linear independence of multi-variable functions in vector space structure.  The beautiful aspect of the method lies in the manipulation of initial approximations of fractional homotopy perturbations towards the approximate solution of the problem, and its accelerated convergence rates. The scheme when applied on model fractional wave equations provided accurate and approximate solutions with ease. \\

\clearpage
\end{document}